\documentclass[oneside,a4paper]{amsart} 

\usepackage[T1]{fontenc}
\usepackage{lmodern}
\usepackage[ngerman,main=english]{babel}  
\usepackage[autostyle]{csquotes}
\usepackage{microtype}
\usepackage{amsmath}
\usepackage{mathtools}
\usepackage{amssymb}
\usepackage[mathscr]{euscript}            
\usepackage{mathrsfs}                     
\usepackage{colonequals}
\usepackage{centernot}
\usepackage{wasysym}
\usepackage[only,llbracket,rrbracket,mapsfrom,longmapsfrom]{stmaryrd}
\usepackage{graphicx}
\usepackage{tikz-cd}
\usetikzlibrary{babel}                    
\usepackage{array}
\usepackage{longtable}
\usepackage{booktabs}
\usepackage{pdfpages}
\usepackage{appendix}
\usepackage{enumitem}
\usepackage[a4paper,textwidth=15cm,textheight=22cm,top=3cm,headheight=15.2pt]{geometry}
\ifdefined\UseBibTeX\else
  \usepackage[backend=biber,sorting=nyt,maxbibnames=99]{biblatex}
\fi
\usepackage{thm-restate}                  
\usepackage[hypertexnames=false]{hyperref}                     
\usepackage[capitalise,nameinlink]{cleveref}

  \defbibheading{bibliography}[\refname]{\section*{#1}} 
\theoremstyle{plain}

\newtheorem{prop}{Proposition}[section]
\newtheorem{lemma}[prop]{Lemma}
\newtheorem{cor}[prop]{Corollary}
\newtheorem{fct}[prop]{Fact}

\theoremstyle{definition}
\newtheorem{definition}[prop]{Definition}

\theoremstyle{remark}
\newtheorem{rem}[prop]{Remark}

\theoremstyle{plain}

\theoremstyle{definition}

\theoremstyle{remark}
\newtheorem*{Rem}{Remark}

\crefname{thm}{Theorem}{Theorems}
\crefname{lemma}{Lemma}{Lemmas}
\crefname{prop}{Proposition}{Propositions}
\crefname{cor}{Corollary}{Corollaries}
\crefname{fct}{Fact}{Facts}
\crefname{conj}{Conjecture}{Conjectures}
\crefname{question}{Question}{Questions}
\crefname{definition}{Definition}{Definitions}
\crefname{notation}{Notation}{Notation}
\crefname{Exc}{Exercise}{Exercises}
\crefname{rem}{Remark}{Remarks}
\crefname{example}{Example}{Examples}

  {\begin{proof}}{\end{proof}}

  {\begin{proof}[Beweisidee]}{\end{proof}}
  
  {\par\medskip\noindent
   \textit{Claim\if\relax\detokenize{#1}\relax\else\ #1\fi.}\ \ignorespaces}%
  {\par}
  {\begin{enumerate}[label=(\alph*),topsep=0pt,partopsep=10pt,itemsep=1ex]}%
  {\end{enumerate}}
\newenvironment{myenumerateromnum}%
  {\begin{enumerate}[label=(\roman*),topsep=0pt,partopsep=5pt,itemsep=0.5ex]}%
  {\end{enumerate}}

\let\phi\varphi
\let\theta\vartheta
\let\epsilon\varepsilon

\let\leq\leqslant
\let\geq\geqslant

\NewCommandCopy{\textOslash}{\O}
\NewCommandCopy{\textSS}{\SS}
\DeclareRobustCommand{\O}{\TextOrMath{\textOslash}{\mathrm{O}}}
\DeclareRobustCommand{\SS}{\TextOrMath{\textSS}{\scriptstyle}}

\newcommand{\actson}{\curvearrowright}

 \newcommand{\BN}{\mathbb{N}}

 \newcommand{\CB}{\mathcal{B}} \newcommand{\CC}{\mathcal{C}} \newcommand{\CD}{\mathcal{D}}

   \newcommand{\CP}{\mathcal{P}}

\newcommand{\FA}{\mathfrak{A}}

\newcommand{\Orb}{\mathcal{O}}

\DeclareMathOperator{\Aut}{Aut}

\DeclareMathOperator{\Pol}{Pol}

\DeclareMathOperator{\CSP}{CSP}

\DeclarePairedDelimiterX{\setof}[2]{\{}{\}}{#1 \;\delimsize\vert\; #2}

\newcommand{\bigtimes}{\mathop{\raisebox{-2pt}{\huge$\times$\kern-1pt}}}
\newcommand{\bigboxplus}{\mathop{\raisebox{-2pt}{\scalebox{1.6}{$\kern1pt\boxplus\kern1pt$}}}}

\newcommand{\beh}{\text{beh}}

\newcommand{\eqG}{\approx_G}

\title{Automatic Continuity for Projective Clone Homomorphisms}

\author{Leopold Karl}
\author{Michael Pinsker}

\usepackage{orcidlink}
\thanks{\hspace*{-5.2mm} Both authors: \textbf{Institut f\"ur Diskrete Mathematik \& Geometrie, Technische Universit\"at Wien}, Vienna, Austria.\\
E-mail addresses: \texttt{leopold.karl@tuwien.ac.at}, \texttt{michael.pinsker@tuwien.ac.at}.\\
ORCID iDs: 0009-0002-0083-4868 \orcidlink{0009-0002-0083-4868} (Leopold Karl), 0000-0002-4727-918X \orcidlink{0000-0002-4727-918X} (Michael Pinsker).\vspace*{1mm}\\
Funded by the European Union (ERC, \textbf{POCOCOP}, 101071674). Views and opinions expressed are however those
of the author(s) only and do not necessarily reflect those of the European Union or the European Research Council
Executive Agency. Neither the European Union nor the granting authority can be held responsible for them. This
research was funded in whole or in part by the Austrian Science Fund (FWF) [I 5948]. For the purpose of Open
Access, the authors have applied a CC BY public copyright licence to any Author Accepted Manuscript (AAM)
version arising from this submission.\vspace*{1mm}\\
\textbf{AI disclosure.} Anthropic's Claude was used for proof search and final LaTeX editing;
it contributed to the idea of Lemma~\ref{lem:cognatedtosimilarcanonical}. The authors take full responsibility for the paper.}

\begin{document}

\begin{abstract}
  We show that every projective clone homomorphism from a closed function clone containing the automorphism group of a countable homogeneous Ramsey structure $\Delta$ in a finite relational language is uniformly continuous. This answers, for such structures,  open questions previously raised by Barto, Bodirsky, Marimon, Pongr\'acz, and the second author in several papers (see \cite{BodirskyPinskerPongracz2021}, \cite{MarimonPinsker2025}, \cite{BartoPinsker2020}). Consequently, the polymorphism clone of a first-order reduct $\Gamma$ of $\Delta$ admits a projective clone homomorphism if and only if $\Gamma$ primitively positively interprets all finite structures; in other words, this ability of $\Gamma$ is characterised by the absence of non-trivial identities in its polymorphism clone, i.e., a purely algebraic condition.
\end{abstract}

\maketitle

\section{Introduction}
\label{sec:intro}

\subsection{Primitive positive interpretations and CSPs}

A \emph{constraint satisfaction problem} asks whether a finite collection of local requirements can be satisfied simultaneously. Formally, one fixes a relational structure $\Gamma$ in a finite relational language, called the \emph{template}; the problem $\CSP(\Gamma)$ then takes as input a finite structure in the same language and asks whether it admits a homomorphism to $\Gamma$. Equivalently, the input is a conjunction of atomic formulas, and the task is to decide whether it is satisfiable in $\Gamma$. For example, the complete graph $K_3$ as a template yields graph $3$-colourability. Varying the template changes the available constraints and, with them, the complexity of the problem.

To compare templates, one needs a way of expressing the constraints of one within another, and \emph{primitive positive interpretations} provide such a mechanism. A \emph{primitive positive} (pp) formula is an existentially quantified conjunction of atomic formulas, and a pp-interpretation represents one structure inside a finite power of another: its domain is a quotient of a pp-definable set, and both equality on this quotient and its relations have pp-definable preimages. If $\Gamma$ pp-interprets $\Lambda$, then $\CSP(\Lambda)$ reduces to $\CSP(\Gamma)$ in polynomial time, simply by replacing variables with tuples and constraints with their defining formulas; this principle is fundamental to the algebraic approach of Bulatov, Jeavons, and Krokhin~\cite{BJK2005}. A template that pp-interprets all finite structures can therefore express every finite-domain CSP, and its own CSP is NP-hard.

For finite templates, this expressive power marks precisely the hard side of the dichotomy. Feder and Vardi conjectured that every finite-domain CSP is in P or NP-complete~\cite{FV1998}, and the algebraic approach ties this distinction to the \emph{polymorphisms} of the template: the operations preserving all of its relations or, equivalently, the homomorphisms $\Gamma^k \to \Gamma$ for $k \geq 1$. The polymorphisms form a \emph{clone} $\Pol(\Gamma)$, that is, a set of operations which contains the projections and is closed under composition. After passing to a core and adding all unary singleton relations, which preserves the CSP up to polynomial-time equivalence, one may assume that this clone is idempotent, meaning that $f(x,\ldots,x) = x$ for every polymorphism $f$~\cite{BJK2005}. In this setting, the algebraic dichotomy conjecture of Bulatov, Jeavons, and Krokhin predicted NP-completeness if the template pp-interprets all finite structures, and polynomial-time solvability otherwise~\cite{BJK2005}. The independent proofs of Bulatov and Zhuk confirmed this prediction and, with it, the Feder--Vardi conjecture~\cite{Bulatov2017,Zhuk2017}.

Crucially, the interpretation condition has a purely algebraic description. Write $\CP$ for the clone of projections on a two-element set. A finite structure $\Gamma$ pp-interprets all finite structures if and only if there is a clone homomorphism $\Pol(\Gamma) \to \CP$, that is, a map preserving arities, projections, and composition~\cite{BJK2005,BodirskyPinsker2015}. Such a homomorphism is called \emph{projective}. Its existence means that the polymorphisms of $\Gamma$ satisfy no non-trivial identities: every finite system of identities they satisfy is also satisfied by projections. The expressive power of a finite template is thus visible in the identities of its polymorphism clone alone.

\subsection{Infinite templates and the role of topology}

Much of this programme remains possible for infinite templates, provided they are countable and $\omega$-categorical. The automorphism group of such a structure has only finitely many orbits on tuples of each fixed length, which provides a form of finiteness despite the infinite domain. The corresponding dichotomy conjecture, due to Bodirsky and the second author, concerns first-order reducts with finite signature of finitely bounded homogeneous structures. Here, a structure is \emph{homogeneous} if every isomorphism between its finite substructures extends to an automorphism, and \emph{finitely bounded} if its finite substructures are characterised by finitely many forbidden induced substructures. The conjecture predicts NP-completeness if the model-complete core of the template, expanded by finitely many constants, pp-interprets all finite structures, and membership in P otherwise~\cite[Conjecture~1.2]{BodirskyPinskerPongracz2021}.

As in the finite case, the interpretation condition can be characterised in terms of polymorphisms, and in fact in much greater generality. For every countable $\omega$-categorical structure $\Gamma$, \emph{Topological Birkhoff}~\cite{BodirskyPinsker2015} shows that $\Gamma$ pp-interprets all finite structures if and only if $\Pol(\Gamma)$ admits a \emph{continuous} projective clone homomorphism. The topology is that of pointwise convergence on operations, with the domain discrete; continuity thus requires the image of an operation to be determined by finitely many of its values. Over finite domains this requirement is automatic, but over infinite domains it is an additional condition, and one that is not algebraic in nature. This raises the following question.

\begin{quote}
    \emph{Does the existence of a projective clone homomorphism from $\Pol(\Gamma)$ imply the existence of a continuous one?}
\end{quote}

For general countable $\omega$-categorical structures $\Gamma$, this question is open~\cite{BodirskyPinskerPongracz2021,MarimonPinsker2025}. A positive answer would recover, for parameter-free pp-interpretability, the purely algebraic criterion that holds over finite domains.

\subsection{Variants of the question}

The importance of this criterion has led to several closely related questions. Their answers are instructive, and they show why the precise form of the question matters.

First, one can allow parameters in interpretations. For countable $\omega$-categorical model-complete cores, Barto and the second author proved in \emph{Topology is irrelevant} that a purely algebraic characterisation is then available~\cite{BartoPinsker2020}. More precisely, the polymorphism clone of some expansion by finitely many constants admits a projective clone homomorphism if and only if the polymorphism clone of some such expansion admits a continuous one; the two tuples of constants need not coincide. Equivalently, the structure pp-interprets all finite structures with parameters if and only if its polymorphism clone has no pseudo-Siggers operation, that is, no operations $s$ of arity six and $\alpha, \beta$ of arity one satisfying
\[
\alpha s(x,y,x,z,y,z) = \beta s(y,x,z,x,z,y).
\]
Naming constants thus yields a purely algebraic criterion, but does not settle the question for interpretations without parameters.

Second, one can enlarge the class of reductions to \emph{pp-constructions}, which additionally permit passing to homomorphically equivalent structures. For countable $\omega$-categorical structures, pp-constructibility of all finite structures corresponds to the existence of a uniformly continuous \emph{minion homomorphism} to $\CP$~\cite{BartoOprsalPinsker2018}. Such a map preserves arities and minors, i.e.\ permuting variables, identifying variables, and adding dummy variables, but need not preserve arbitrary composition. In this setting topology cannot be omitted in general: \emph{Topology is relevant} constructs an $\omega$-categorical structure whose polymorphism clone admits a minion homomorphism to $\CP$, but no uniformly continuous one~\cite{BMOOPW2019}. The example, however, lies outside the class of reducts of finitely bounded homogeneous structures, and hence outside the scope of the Bodirsky--Pinsker conjecture.

Neither result resolves the question itself, which appears in \cite{BodirskyPinskerPongracz2021} and \cite{BartoPinsker2020}, and again in the recent guide to topological reconstruction~\cite{MarimonPinsker2025}: for projective clone homomorphisms, with no constants added, does existence imply the existence of a continuous one? Its persistence reflects the gap between the algebraic information carried by identities and the finite descriptions required by pp-interpretations.

\subsection{The mathematical appeal of automatic continuity}

The question is of interest well beyond constraint satisfaction. To begin with, the gap just described has a classical universal-algebraic formulation. Let $\FA_\Gamma$ be the algebra on the domain of $\Gamma$ which has one fundamental operation for every member of $\Pol(\Gamma)$. A projective clone homomorphism exists precisely when the variety $\mathrm{HSP}(\FA_\Gamma)$ generated by $\FA_\Gamma$ contains a two-element algebra all of whose fundamental operations are projections~\cite{Birkhoff-abstract-algebras,BodirskyPinskerPongracz2021}. For countable $\omega$-categorical $\Gamma$, a continuous projective clone homomorphism exists precisely when such an algebra already belongs to the pseudovariety $\mathrm{HSP}_{\mathrm{fin}}(\FA_\Gamma)$, which consists of the homomorphic images of subalgebras of finite powers of $\FA_\Gamma$~\cite{BodirskyPinsker2015}. Our question thus asks whether a two-element projection algebra in the variety forces one in the pseudovariety. Viewed this way, it concerns the passage from arbitrary powers to finite powers, and connects CSP theory to questions about algebraic representation that long predate it.

Moreover, the question belongs to a broader theme: the reconstruction of topological from algebraic information. By the theorem of Ahlbrandt and Ziegler, the automorphism group of a countable $\omega$-categorical structure, viewed as a topological group, determines the structure up to first-order bi-interpretability~\cite{AZ1986}, and topological polymorphism clones play the corresponding role for pp-bi-interpretability~\cite{BodirskyPinsker2015}. This motivates the study of when the abstract algebraic structure of a group or clone already determines the relevant topological information; see~\cite{BPP2017,MarimonPinsker2025}.

From this perspective, the strongest conceivable answer to our question is \emph{automatic continuity with respect to $\CP$}: every projective clone homomorphism is continuous. This asks for more than the existence of a continuous projective clone homomorphism whenever an abstract one exists, and for general $\omega$-categorical structures it is too much to ask. Indeed, Bodirsky, Pongr\'acz, and the second author exhibited an $\omega$-categorical structure, first-order definable in a Cherlin--Hrushovski structure, whose polymorphism clone has a discontinuous projective clone homomorphism~\cite{BodirskyPinskerPongracz2021}.

Their construction has a group-theoretic counterpart. The automorphism group of the Cherlin--Hrushovski structure has a continuous quotient $(\mathbb Z/2\mathbb Z)^{\mathbb N}$~\cite{Tsankov2012}, and a non-principal ultrafilter yields a discontinuous homomorphism from this product onto $\mathbb Z/2\mathbb Z$. At the level of clones, an ultrafilter similarly selects argument positions of operations on two-element quotients. This is a common construction principle rather than an extension of the group homomorphism, since a projective clone homomorphism necessarily sends every unary operation to the identity. The counterexample does not rule out automatic continuity for templates within the scope of the Bodirsky--Pinsker conjecture; it does show, however, that any positive result must use more than $\omega$-categoricity.

\subsection{Our contribution}

We answer the question in a strong form for a large class of structures: we establish automatic continuity whenever the template admits a homogeneous Ramsey expansion in a finite relational language. Notably, we do not require finite boundedness.

\begin{restatable}{thm}{MainThm}\label{thm:main}
  Let $\Delta$ be a countable homogeneous Ramsey structure in a finite relational language and $\CD$ a closed function clone containing $\Aut(\Delta)$. Then every projective clone homomorphism $\xi\colon \CD \to \CP$ is uniformly continuous.\newline
  In particular, the polymorphism clone of every first-order reduct of $\Delta$ has automatic continuity with respect to $\CP$.
\end{restatable}

For closed function clones containing $\Aut(\Delta)$, Theorem~\ref{thm:main} answers \cite[Questions~1.3 and~2.1]{BodirskyPinskerPongracz2021}, the question in \cite[Sec.~6.1]{BartoPinsker2020} and the first part of \cite[Question~7.2]{MarimonPinsker2025} positively. Moreover, it shows that for structures first-order definable in such a $\Delta$ there are no discontinuous projective clone homomorphisms, which concerns the last part of \cite[Question~7.2]{MarimonPinsker2025}.

Combined with known results, Theorem~\ref{thm:main} yields the following corollary. The equivalence of the primed items follows from Birkhoff's HSP theorem~\cite{Birkhoff-abstract-algebras}, and that of the non-primed items was shown in~\cite{BodirskyPinsker2015}; Theorem~\ref{thm:main} connects the two groups.

\begin{restatable}{cor}{MainCor}\label{cor:main}
    Let $\Gamma$ be a first-order reduct of a structure $\Delta$ as in Theorem~\ref{thm:main}, and let $\FA_\Gamma$ be the algebra on the domain of $\Gamma$ whose operations are the members of $\Pol(\Gamma)$. Then the following are equivalent.
    \begin{myenumerateromnum}
        \item[(i)] $\Pol(\Gamma)$ has a uniformly continuous projective clone homomorphism.
        \item[(i')] $\Pol(\Gamma)$ has a projective clone homomorphism.
        \item[(ii)] The pseudovariety generated by $\FA_\Gamma$ contains a two-element algebra all of whose operations are projections.
        \item[(ii')] The variety generated by $\FA_\Gamma$ contains a two-element algebra all of whose operations are projections.
        \item[(iii)] $\Gamma$ primitively positively interprets all finite structures.
    \end{myenumerateromnum}
\end{restatable}

In particular, for all structures within the scope of the Bodirsky--Pinsker conjecture that admit a homogeneous Ramsey expansion in a finite relational language, Corollary~\ref{cor:main} provides a purely algebraic characterisation of the ability to pp-interpret all finite structures.

The Ramsey assumption is natural in this context, as it connects our theorem to an established approach to infinite-domain CSPs. Ramsey expansions allow arbitrary polymorphisms to be studied through \emph{canonical} operations, whose behaviour is controlled by the orbits of tuples. Besides their role in the analysis of CSP templates, these methods underlie the decidability results for pp-definability in~\cite{BPT-decidability-definability} and for pp-bi-interpretability in~\cite{FP2026}. They also interact well with the passage to cores: Mottet and the second author proved that the existence of an $\omega$-categorical homogeneous Ramsey expansion is preserved by taking the model-complete core, as is finite boundedness when present~\cite{MottetPinsker-cores-ramsey}. Accordingly, the Ramsey property is a standard additional assumption for structures within the scope of the Bodirsky--Pinsker conjecture, and the existing complexity classifications concern reducts of structures which have it. In fact, it is an open question whether every countable homogeneous structure in a finite relational language admits a homogeneous Ramsey expansion in a finite relational language~\cite{BPT-decidability-definability,Bodirsky-ramsey-classes,NguyenVanThe-KPT-correspondence,MottetPinsker-cores-ramsey}. If the answer is positive, then our result yields automatic continuity for all structures within the scope of the conjecture, and beyond.

The key notion of our proof is \emph{cognacy}. Call two operations \emph{locally equal modulo $\Aut(\Delta)$} if on every finite set they agree up to composition with an automorphism of $\Delta$ from the outside. Two operations in $\CD$ are cognate if, up to this equivalence, they arise from a common operation in $\CD$ by composing its arguments with elements of $\overline{\Aut(\Delta)}$, the closure of $\Aut(\Delta)$ in the space of all unary functions. Since a projective clone homomorphism sends every unary operation to the identity, it cannot distinguish cognate operations. Using a compactness argument and the canonisation theorem for extremely amenable groups~\cite{BodirskyPinsker2021}, we show that every operation in $\CD$ is cognate with a canonical one. As there are, in each arity, only finitely many canonical operations up to local equality modulo $\Aut(\Delta)$, cognacy with canonical operations is already decided on a finite set. This yields continuity of the homomorphism, and uniform continuity then follows by a known compactness argument.

\subsection{Organisation of the paper}
Section~\ref{sec:prelim} recalls the necessary notions. Section~\ref{sec:equalitymoduloG} recalls local equality modulo a permutation group, following \cite{BodirskyPinsker2021}, and Section~\ref{sec:cognacy} introduces the notion of cognacy. In Section~\ref{sec:ramsey}, we show that every operation is cognate with a canonical one, using the Ramsey property of $\Delta$. The Main Theorem~\ref{thm:main} and its Corollary~\ref{cor:main} are proved in Section~\ref{sec:automaticcontinuity}.

\section{Preliminaries}
\label{sec:prelim}

We only recall the notions used in this paper and refer to \cite{Hodges1997} for model theory, to \cite{BodirskyPinskerPongracz2021,BartoOprsalPinsker2018} for function clones and their topology, and to \cite{BourbakiTopology1} for uniform spaces.

\emph{Permutation groups and structures.}
Any permutation group $G \actson D$ acts componentwise on $D^m$, and $G\actson D$ is \emph{oligomorphic} if this action has only finitely many orbits for every $m \geq 1$. We write $\overline{G} \subseteq D^D$ for the closure of $G$ in the space of all functions of $D^D$ (not just bijections) with respect to the topology of pointwise convergence; $\overline{G}$ is a monoid of injections.\newline
A relational structure $\Delta$ is \emph{homogeneous} if every isomorphism between finite substructures of $\Delta$ extends to an automorphism of $\Delta$. If $\Delta$ is homogeneous in a finite relational language, then $\Aut(\Delta)$ is oligomorphic, since the orbit of a tuple is determined by its quantifier-free type; hence $\Delta$ is $\omega$-categorical by the theorem of Engeler, Ryll-Nardzewski and Svenonius \cite{Hodges1997}. A \emph{first-order reduct} of $\Delta$ is a structure on the same domain all of whose relations are first-order definable in $\Delta$ without parameters; its automorphism group contains $\Aut(\Delta)$.\newline
The \emph{age} of $\Delta$ is the class of finite structures that embed into $\Delta$. It has the \emph{Ramsey property} if for all $A, B$ in the age there is some $C$ in the age such that for every colouring of the embeddings $A \to C$ with two colours there exists an embedding $e\colon B \to C$ such that all embeddings $e \circ a$, where $a\colon A \to B$ is an embedding, have the same colour. We call $\Delta$ \emph{Ramsey} if its age has the Ramsey property. By the Kechris--Pestov--Todor\v{c}evi\'c correspondence \cite{KPT05}, a countable homogeneous structure is Ramsey if and only if its automorphism group is extremely amenable; see also \cite[Sec.~1]{BodirskyPinsker2021}.

\emph{Function clones and topology.}
A \emph{function clone} $\CC$ on $D$ is a set of finitary operations on $D$ which is closed under composition and contains all projections $\pi_i^k$; we write $\CC^{(k)}$ for its $k$-ary part. The \emph{polymorphism clone} $\Pol(\Gamma)$ of a structure $\Gamma$ consists of all finitary operations preserving all relations of $\Gamma$. We equip $D^{D^k}$ with the topology of pointwise convergence, i.e.\ the product topology where $D$ is discrete, and the set of all finitary operations on $D$ with the corresponding sum topology. A function clone is \emph{closed} if it is closed in this space; the closed function clones are precisely the polymorphism clones \cite[Cor.~6.1.6]{Bodirsky2021}. In particular, if $\Gamma$ is a first-order reduct of $\Delta$, then $\Pol(\Gamma)$ is a closed function clone containing $\Aut(\Delta)$.\newline
Analogously, $D$ carries the discrete uniformity, $D^{D^k}$ the corresponding product uniformity, whose basic entourages are the sets $\{(f,g) \mid f|_F = g|_F\}$ for finite $F \subseteq D^k$, and the set of all finitary operations the corresponding sum uniformity.

\emph{Clone homomorphisms.}
A \emph{clone homomorphism} $\xi\colon \CC \to \CB$ between function clones preserves arities, sends each projection $\pi_i^k$ to $\pi_i^k$, and preserves composition. We write $\CP$ for the clone of projections on $\{0,1\}$, which is finite in each arity and discrete. Clone homomorphisms to $\CP$ are called \emph{projective}; every projective clone homomorphism sends all unary operations to $\pi_1^1$. Since $\CP$ is discrete, an arity-preserving map $\xi\colon \CC \to \CP$ is
\begin{itemize}
    \item \emph{continuous} if and only if for every $f \in \CC^{(k)}$ there is a finite $F \subseteq D^k$ such that $\xi(g) = \xi(f)$ for all $g \in \CC^{(k)}$ with $g|_F = f|_F$, and
    \item \emph{uniformly continuous} if and only if for every $k \geq 1$ there is a finite $F \subseteq D^k$ such that $\xi(f) = \xi(g)$ for all $f,g \in \CC^{(k)}$ with $f|_F = g|_F$.
\end{itemize}
If $D$ is finite, every such map is uniformly continuous; hence we may and do assume that $D$ is countably infinite.

\begin{fct}[{\cite{BodirskyPinsker2015}, see \cite[Prop.~1.1]{BodirskyPinskerPongracz2021}}]\label{fct:BP}
    Let $\Gamma$ be a countable $\omega$-categorical structure. Then $\Pol(\Gamma)$ has a continuous projective clone homomorphism if and only if $\Gamma$ primitively positively interprets all finite structures.
\end{fct}

\section{Equality modulo a Permutation Group}
\label{sec:equalitymoduloG}

This section provides a setting that unifies the various compactness arguments throughout the paper. The notion of equality modulo a permutation group goes back to \cite{BodirskyPinsker2015}, 
and we present it in the form of \cite[Sec.~2]{BodirskyPinsker2021}.

\begin{definition}
    \label{def:locEqModG}
    Let $G \actson Y$ be a permutation group, and $f,g \in Y^X$ for some set $X$. We say that $f$ locally equals $g$ modulo $G$, and write $f \eqG g$, if for every finite $F \subseteq X$ there exist $\alpha_1,\alpha_2 \in G$ such that $\alpha_1 \circ f|_F = \alpha_2 \circ g|_F$.\newline
    We say that $f$ equals $g$ globally modulo $\overline{G}$, where $\overline{G} \subseteq Y^Y$ denotes the closure of $G$ in the topology of pointwise convergence, if there exist $e_1, e_2 \in \overline{G}$ with $e_1 \circ f = e_2 \circ g$.
\end{definition}

\begin{Rem}
    Since every element of $\overline{G}$ agrees with some element of $G$ on any given finite set, global equality modulo $\overline{G}$ implies local equality modulo $G$. Conversely, by Lemma~\ref{lem:lift} below, local equality modulo $G$ implies global equality modulo $\overline{G}$ whenever $G$ is oligomorphic and $X, Y$ are countable.\newline
    Also, notice that local equality modulo $G$ is an equivalence relation, since $G$ is a group. For the same reason, $f \eqG g$ holds if and only if for every finite $F \subseteq X$ there exists some $\alpha \in G$ with $\alpha \circ f|_F = g|_F$, i.e.\ if and only if $g$ lies in the closure of $\{\alpha \circ f \mid \alpha \in G\}$. We will mostly use this reformulation.
\end{Rem}

Throughout, $Y$ carries the discrete topology, $Y^X$ the topology of pointwise convergence (i.e.\ the product topology), and $Y^X/_{\eqG}$ the corresponding quotient topology. For $f \in Y^X$ we write $[f] = [f]_{\eqG}$ for its equivalence class and $\pi_{\eqG}\colon Y^X \to Y^X/_{\eqG}$, $f \mapsto [f]$, for the quotient map. For finite $F \subseteq X$ we write $V_F(f) = \{g \in Y^X \mid g|_F = f|_F\}$; these sets form a base of the topology of $Y^X$. Notice that $f \eqG g$ implies $f|_F \eqG g|_F$ in $Y^F$ for every $F \subseteq X$.

The following lemma is the main motivation for introducing the notion of local equality modulo $G$.

\begin{lemma}[{\cite[Lem.~4]{BodirskyPinsker2021}}]
    \label{lem:compactness}
    If $G \actson Y$ is oligomorphic, and $X$ is countable, then the space $Y^X/_{\eqG}$ is compact and Hausdorff.
\end{lemma}

A proof can be found in \cite[Lem.~4]{BodirskyPinsker2021}. The idea is to show that the map $[f] \mapsto \bigl([f|_F]\bigr)_{F \in [X]^{<\omega}}$ is a homeomorphism of $Y^X/_{\eqG}$ onto a closed subspace of the compact Hausdorff space $\prod_{F \in [X]^{<\omega}} Y^F/_{\eqG}$, whose factors are finite and discrete as $G$ is oligomorphic.

Local equality modulo $G$ only expresses agreement on finite sets. The following lemma is a consequence of the compactness of the spaces $Y^X/_{\eqG}$ (Lemma~\ref{lem:compactness}).

\begin{lemma}[{\cite[Lem.~3]{BodirskyPinsker2021}}]
    \label{lem:lift}
    Let $G \actson Y$ be oligomorphic with $Y$ countable, let $X$ be countable and let $f,g \in Y^X$ with $f \eqG g$. Then $f$ equals $g$ globally modulo $\overline{G}$.
\end{lemma}

Finally, the following lemma allows us to lift convergence of a sequence in $Y^X/_{\eqG}$ to the convergence of a sequence in $Y^X$. A proof (in the special case of Siggers operations) is found in \cite[Lemma~4.2]{BartoPinsker-algebraic-dichotomy}.

\begin{lemma}
\label{lem:liftconvergence}
    Let $G \actson Y$ be a permutation group, $X$ countable, $(f_n)_{n \in \BN}$ a sequence in $Y^X$ and $f \in Y^X$. Then the following are equivalent:
    \begin{myenumerateromnum}
        \item The sequence $([f_n])_{n \in \BN}$ converges to $[f]$ in $Y^X/_{\eqG}$.
        \item There exists a sequence $(\alpha_n)_{n \in \BN}$ in $G$ such that $(\alpha_n \circ f_n)_{n \in \BN}$ converges to $f$ in $Y^X$.
    \end{myenumerateromnum}
\end{lemma}



\section{Cognacy}
\label{sec:cognacy}

The goal of this section is to identify sets of operations on which every projective clone homomorphism $\xi$ is necessarily constant. To this end, we exploit a property specific to projective clone homomorphisms: they map all unary operations to the same element $\pi_1^1$ of the clone of projections $\CP$. The notion of cognacy is built on this property and yields a sufficient condition for the constancy of projective clone homomorphisms on certain sets: whenever two operations are cognate, their images under $\xi$ agree (see Corollary~\ref{cor:xi}).

Let $D$ be a countably infinite set, $G \actson D$ an oligomorphic permutation group and $\overline{G} \subseteq D^D$ its closure in the topology of pointwise convergence. Let $\CD$ be a closed function clone on $D$ that contains $G$. Since $\CD$ is closed, it also contains $\overline{G}$. We fix an enumeration $d_0, d_1, \dots$ of $D$ and set $D_n = \{d_0, \dots, d_n\}$.

\begin{definition}\label{def:linked}
The operations $f,g \in \CD^{(k)}$ are \emph{cognate} if there are $L\in\CD^{(k)}$ and elements $e_1,\dots,e_k,e_1',\dots,e_k'\in\overline{G}$ such that
\begin{equation}\label{eq:linked}
f \eqG L(e_1,\dots, e_k)\qquad\text{and}\qquad L(e_1',\dots,e_k') \eqG g.
\end{equation}
We call $L$ a \emph{cognacy witness} for $f$ and $g$.\newline
For a finite set $F \subseteq D$ the operations $f$ and $g$ are \emph{F-cognate} if there is $L_F \in \CD^{(k)}$ such that both identities in \eqref{eq:linked} hold on $F^k$.\newline
The operations $f$ and $g$ are \emph{finitely cognate} if they are $F$-cognate for every finite $F \subseteq D$.
\end{definition}

\begin{rem}\label{rem:eqGgivescognacy}
    Notice that cognacy is symmetric, and that $f \eqG g$ implies that $f$ and $g$ are cognate, as witnessed by $L = f$ and $e_i = e'_i = \mathrm{id}$.
\end{rem}

The following lemma is an immediate consequence of the above definitions.

\begin{lemma}
    \label{lem:eqModGandLinkedness}
    If $f \eqG g$ and $g$ is (finitely) cognate with $h$, then $f$ is (finitely) cognate with $h$.
\end{lemma}

\begin{proof}
    The same witness establishes cognacy, and the same family of witnesses establishes finite cognacy.
\end{proof}

As a corollary, we obtain that the notion of (finite) cognacy is well-defined on the equivalence classes of local equality modulo $G$ from Section~\ref{sec:equalitymoduloG}.

\begin{cor}
    For $f,g \in \CD^{(k)}$ the following are equivalent.
    \begin{myenumerateromnum}
        \item $f$ and $g$ are (finitely) cognate.
        \item For all $f' \in [f]$ and $g' \in [g]$ the operations $f'$ and $g'$ are (finitely) cognate.
    \end{myenumerateromnum}
    Whenever these conditions are satisfied, we say that $[f]$ and $[g]$ are (finitely) cognate.
\end{cor}

\begin{proof}
    The implication $(ii) \implies (i)$ is immediate. Suppose $f$ and $g$ are (finitely) cognate and let $f' \in [f]$ and $g' \in [g]$. By applying Lemma~\ref{lem:eqModGandLinkedness} twice, using the symmetry of (finite) cognacy, we obtain that $f'$ and $g'$ are (finitely) cognate.
\end{proof}

We now use the setting of Section~\ref{sec:equalitymoduloG} to show that the notions of finite cognacy and cognacy agree.

\begin{prop}\label{prop:universalcognacy}
Let $f,g \in \CD^{(k)}$. Then $f$ and $g$ are finitely cognate if and only if they are cognate.
\end{prop}

\begin{proof}
    Clearly, cognacy implies finite cognacy.\newline
    Conversely, let $f$ and $g$ be finitely cognate and let $L_n \in \CD^{(k)}$ and $\alpha_{n,j}, \alpha'_{n,j} \in G$ be finite cognacy witnesses for $F = D_n$, i.e.\
\begin{equation}\label{eq:finLinked}
    f \eqG L_n(\alpha_{n,1},\dots, \alpha_{n,k})\qquad\text{and}\qquad L_n(\alpha'_{n,1}, \dots,\alpha'_{n,k}) \eqG g
\end{equation}
    hold on $D_n^k$. 

\emph{Inner normalisation.}
    We first normalise all $\alpha_{n,j}$ and $\alpha'_{n,j}$ simultaneously by a single element of $G$ from the outside. To this end, let $G$ act diagonally on $D^{2k}$, which is again oligomorphic, and consider the functions
    $$a_n\colon D \to D^{2k},\quad x \mapsto \left(\alpha_{n,1}(x), \dots, \alpha_{n,k}(x), \alpha'_{n,1}(x), \dots, \alpha'_{n,k}(x)\right).$$
    By Lemma~\ref{lem:compactness}, the sequence $([a_n])_{n \in \BN}$ in $(D^{2k})^D/_{\eqG}$ admits a subsequence converging to some $[e]$. Applying the convergence lifting Lemma~\ref{lem:liftconvergence} yields a sequence $(\sigma_n)_{n \in \BN}$ in $G$ with $\sigma_n \circ a_n$ converging to $e$ in $(D^{2k})^D$ along the subsequence. Writing $e = (e_1, \dots, e_k, e'_1, \dots, e'_k)$, this means that $\sigma_n \circ \alpha_{n,j} \to e_j$ and $\sigma_n \circ \alpha'_{n,j} \to e'_j$ along the subsequence for all $j \leq k$. Since they are limits of elements of $G$, we have $e_j, e'_j \in \overline{G}$.

\emph{A global witness.}
    After reindexing the convergent subsequence, define $L_n^* = L_n(\sigma_n^{-1}, \dots, \sigma_n^{-1}) \in \CD^{(k)}$. These functions are finite cognacy witnesses for $D_n$ since they satisfy:

\begin{equation}\label{eq:newfinLinks}
    f \eqG L_n^*(\sigma_n \circ \alpha_{n,1}, \dots, \sigma_n \circ \alpha_{n,k})\qquad\text{and}\qquad L_n^*(\sigma_n \circ\alpha'_{n,1}, \dots,\sigma_n \circ\alpha'_{n,k}) \eqG g
\end{equation}

    on $D_n^k$.
    Again, by the Compactness Lemma~\ref{lem:compactness}, the sequence $([L_n^*])_{n \in \BN}$ in $D^{D^k}/_{\eqG}$ admits a subsequence $[L_i^*]_{i \in I}$ converging to some $[L]$. We claim that $L$ is a cognacy witness for $f$ and $g$ satisfying the identities
\begin{equation}\label{eq:globalLink}
    f \eqG L(e_1,\dots,e_k)\qquad\text{and}\qquad L(e'_1, \dots, e'_k) \eqG g.
\end{equation}

    Indeed, by \eqref{eq:newfinLinks}, the sequence $[L_n^*(\sigma_n \circ \alpha_{n,1}, \dots, \sigma_n \circ \alpha_{n,k})]$ converges to $[f]$ and along the subsequence $([L_i^*])_{i \in I}$ it converges to $[L(e_1,\dots,e_k)]$. Since $D^{D^k}/_{\eqG}$ is Hausdorff by Lemma~\ref{lem:compactness}, we conclude that $[f] = [L(e_1,\dots,e_k)]$, i.e. $f \eqG L(e_1, \dots, e_k)$. The same argument applied to $g$ yields $[g] = [L(e'_1,\dots,e'_k)]$.
\end{proof}

As a corollary, we obtain that even finite cognacy of two operations is sufficient for their images to agree under every projective clone homomorphism.

\begin{cor}\label{cor:xi}
Let $\xi\colon\CD\to\CP$ be a projective clone homomorphism, and let $f$ and $g$ be finitely cognate. Then $\xi(f)=\xi(g)$.
\end{cor}

\begin{proof}
By \cref{prop:universalcognacy}, $f$ and $g$ are cognate. Let $L$ be a cognacy witness. By Lemma~\ref{lem:lift}, $f$ and $L(e_1,\dots,e_k)$ are globally equal modulo $\overline{G}$, and since $\xi$ maps unary operations to $\pi^1_1$, we obtain $\xi(f) = \xi(L(e_1,\dots,e_k)) = \xi(L)$; similarly $\xi(g) = \xi(L)$.
\end{proof}

\section{Ramsey: every operation is cognate with a canonical one}\label{sec:ramsey}

Throughout this section, we work in the setting of Section~\ref{sec:cognacy}, where moreover $G = \Aut(\Delta)$ for a countably infinite homogeneous Ramsey structure $\Delta$ in a finite relational language, all of whose relations have arity at most $r$, for some $r \in \BN_{\geq 2}$. For $s,t \in D^m$ we write $s \equiv t$ if $s$ and $t$ lie in the same orbit of the componentwise action of $G$ on $D^m$, and we write $\Orb_r$ for the set of orbits of the componentwise action of $G$ on $D^r$, which is finite since $G$ is oligomorphic. Operations are applied to tuples componentwise, i.e.\ for $g\colon D^k \to D$ and $s^1, \dots, s^k \in D^m$ we write $g(s^1,\dots,s^k) = \bigl(g(s^1_1,\dots,s^k_1), \dots, g(s^1_m,\dots,s^k_m)\bigr)$. The following notion of canonicity is taken from \cite{BPT-decidability-definability}.

\begin{definition}\label{def:canonical}
An operation $g\colon D^k\to D$ is \emph{canonical} if for all $m\geq1$ and all $s^1,\dots,s^k,t^1,\dots,t^k\in D^m$ with $s^i\equiv t^i$ for all $i \leq k$ we have $g(s^1,\dots,s^k)\equiv g(t^1,\dots,t^k)$. If this holds for all tuples with entries in a set $F \subseteq D$, we say that $g$ is \emph{canonical on $F$}.
\end{definition}

Since $\Delta$ is homogeneous in a finite relational language, the orbit of a tuple is determined by the orbits of its $r$-subtuples. Hence, it suffices to check canonicity for $m = r$, and a canonical $g$ induces a map $\beh(g)\colon \Orb_r^k \to \Orb_r$, its \emph{behaviour}; canonical operations with the same behaviour are locally equal modulo $G$.

Since $\Delta$ is Ramsey, the following is a direct consequence of \cite[Lemma~8]{BPT-decidability-definability}. It is obtained from the Canonisation Theorem \cite[Thm.~5]{BodirskyPinsker2021} using a characterisation of the Ramsey property \cite{KPT05}, as well.

\begin{lemma}\label{lem:Ramsey}
Let $f \in \CD^{(k)}$ and $F \subseteq D$ finite. Then there are $\alpha_1, \dots, \alpha_k \in G$ such that $f(\alpha_1, \dots, \alpha_k)$ is canonical on $F$.
\end{lemma}

\begin{lemma}\label{lem:cognatedtosimilarcanonical}
For every $f \in \CD^{(k)}$ there exists some $N_f \in \BN$ such that a canonical $g \in \CD^{(k)}$ is cognate with $f$ if and only if it is $D_{N_f}$-cognate with $f$. Moreover, if $f' \in \CD^{(k)}$ agrees with $f$ on $D_{N_f}^k$, then some canonical $g \in \CD^{(k)}$ is cognate with both $f$ and $f'$.
\end{lemma}

\begin{proof}
By Lemma~\ref{lem:Ramsey}, for every $n$ there are $\alpha^{(n)}_1, \dots, \alpha^{(n)}_k \in G$ such that $f(\alpha^{(n)}_1, \dots, \alpha^{(n)}_k)$ is canonical on $D_n$. By Lemmas~\ref{lem:compactness} and~\ref{lem:liftconvergence}, a subsequence converges modulo $G$ to some $g \in \CD^{(k)}$, which is canonical; $f$ and $g$ are finitely cognate (with $L_F = f$), hence cognate by Proposition~\ref{prop:universalcognacy}. So every operation in $\CD^{(k)}$ is cognate with a canonical one.

There are only finitely many $\eqG$-classes of canonical operations, and the set of those classes that are $D_n$-cognate with $f$ decreases in $n$. Hence it is constant from some $N_f$ on, and then consists of the classes finitely cognate, i.e.\ by Proposition~\ref{prop:universalcognacy} cognate, with $f$.

For the second statement, let $g$ be canonical and cognate with $f'$. Then $g$ is $D_{N_f}$-cognate with $f'$, hence with $f$, so $g$ is cognate with $f$.
\end{proof}

\section{Automatic Continuity}
\label{sec:automaticcontinuity}

We prove the main theorem and derive Corollary~\ref{cor:main}.

\MainThm*

\begin{proof}
    Let $\xi: \CD \to \CP$ be a clone homomorphism. In \cite[Lem.~3.4]{BodirskyPinskerPongracz2021} it is shown that continuity of a projective clone homomorphism only depends on the continuity of the restriction to the binary fragment $\xi|_{\CD^{(2)}}$. Although the following argument applies to each $k$-ary fragment $\CD^{(k)}$ of the function clone $\CD$, we use this result and restrict attention to $k=2$ for readability.

    \emph{Continuity.} Let $f \in \CD^{(2)}$ and 
    let $f' \in \CD^{(2)}$ agree with $f$ on $D_{N_f}^2$. By Lemma~\ref{lem:cognatedtosimilarcanonical}, the function $f'$ is cognate with some canonical function $g$ which is cognate with $f$. As the images of cognate operations under any projective clone homomorphism agree by Corollary~\ref{cor:xi}, we conclude that $$\xi(f') = \xi(g) = \xi(f).$$ So $\xi$ is constant on the open neighbourhood $V_{D_{N_f}^2}(f) = \bigl\{h \in \CD^{(2)}\mid h|_{D_{N_f}^2} = f|_{D_{N_f}^2}\bigr\}$ of $f$. Hence, $\xi|_{\CD^{(2)}}$ is locally constant, thus continuous, and we conclude that all of $\xi$ is continuous.


    Uniform continuity follows from continuity by a result in \cite{BodirskyPinsker2015} as noted in \cite[Prop.~6.4]{BartoOprsalPinsker2018}.

\end{proof}

\begin{Rem}
    We remark that the result holds verbatim for clone homomorphisms to any idempotent clone (instead of $\CP$).
\end{Rem}

We conclude with the consequences for first-order reducts announced in Corollary~\ref{cor:main}.

\MainCor*

\begin{proof}
    Theorem~\ref{thm:main} proves the equivalence of $(i)$ and $(i')$. The primed equivalences are shown in \cite{Birkhoff-abstract-algebras} and the non-prime equivalences in \cite{BodirskyPinsker2015}.
\end{proof}


\printbibliography

@article{AZ1986,
  author  = {Ahlbrandt, Gisela and Ziegler, Martin},
  title   = {Quasi finitely axiomatizable totally categorical theories},
  journal = {Annals of Pure and Applied Logic},
  volume  = {30},
  pages   = {63--82},
  year    = {1986},
  doi     = {10.1016/0168-0072(86)90037-0}
}

@article{BartoOprsalPinsker2018,author={Barto, Libor and Opr{\v{s}}al, Jakub and Pinsker, Michael},title={The wonderland of reflections},journal={Israel Journal of Mathematics},volume={223},number={1},pages={363--398},year={2018}}

@inproceedings{BartoPinsker-algebraic-dichotomy,
author = {Libor Barto and Michael Pinsker},
title = {The algebraic dichotomy conjecture for infinite domain {C}onstraint {S}atisfaction {P}roblems},
booktitle = {Proceedings of the 31st Annual ACM/IEEE Symposium on Logic in Computer Science (LICS 2016)},
year = {2016},
pages = {615-622},
doi = {10.1145/2933575.2934544},
note = {Preprint arXiv:1602.04353}}

@article{BartoPinsker2020,author={Barto, Libor and Pinsker, Michael},title={Topology is irrelevant (in a dichotomy conjecture for infinite domain constraint satisfaction problems)},journal={SIAM Journal on Computing},volume={49},number={2},pages={365--393},year={2020}}

@article{Birkhoff-abstract-algebras,
author = {Garrett Birkhoff},
title = {On the Structure of Abstract Algebras},
journal = {Mathematical Proceedings of the Cambridge Philosophical Society},
volume = {31},
number = {4},
year = {1935},
pages = {433-454},
doi = {10.1017/S0305004100013463}}

@book{Bodirsky2021,author={Bodirsky, Manuel},title={Complexity of Infinite-Domain Constraint Satisfaction},series={Lecture Notes in Logic},volume={52},publisher={Cambridge University Press},address={Cambridge},year={2021},doi={10.1017/9781107337534}}

@incollection{Bodirsky-ramsey-classes,
author = {Manuel Bodirsky},
title = {{R}amsey classes: examples and constructions},
booktitle = {Surveys in Combinatorics 2015},
series = {London Mathematical Society Lecture Note Series},
volume = {424},
publisher = {Cambridge University Press},
year = {2015},
pages = {1-48},
doi = {10.1017/CBO9781316106853.002},
note = {Preprint arXiv:1502.05146}}

@inproceedings{BMOOPW2019,
  author        = {Bodirsky, Manuel and Mottet, Antoine and
                   Ol{\v{s}}{\'a}k, Miroslav and Opr{\v{s}}al, Jakub and
                   Pinsker, Michael and Willard, Ross},
  title         = {Topology is relevant (in a dichotomy conjecture for
                   infinite-domain constraint satisfaction problems)},
  booktitle     = {34th Annual {ACM/IEEE} Symposium on Logic in
                   Computer Science ({LICS})},
  pages         = {1--12},
  year          = {2019},
  doi           = {10.1109/LICS.2019.8785883},
  eprint        = {1901.04237},
  archivePrefix = {arXiv}
}

@article{BodirskyPinsker2015,author={Bodirsky, Manuel and Pinsker, Michael},title={Topological {B}irkhoff},journal={Transactions of the American Mathematical Society},volume={367},number={4},pages={2527--2549},year={2015}}

@article{BodirskyPinsker2021,author={Bodirsky, Manuel and Pinsker, Michael},title={Canonical functions: a proof via topological dynamics},journal={Contributions to Discrete Mathematics},volume={16},number={2},pages={36--45},year={2021}}

@article{BodirskyPinskerPongracz2021,author={Bodirsky, Manuel and Pinsker, Michael and Pongr{\'a}cz, Andr{\'a}s},title={Projective clone homomorphisms},journal={Journal of Symbolic Logic},volume={86},number={1},pages={148--161},year={2021}}

@article{BPP2017,
  author  = {Bodirsky, Manuel and Pinsker, Michael and
             Pongr{\'a}cz, Andr{\'a}s},
  title   = {Reconstructing the topology of clones},
  journal = {Transactions of the American Mathematical Society},
  volume  = {369},
  number  = {5},
  pages   = {3707--3740},
  year    = {2017},
  doi     = {10.1090/tran/6937}
}

@article{BPT-decidability-definability,
author = {Manuel Bodirsky and Michael Pinsker and Todor Tsankov},
title = {Decidability of definability},
journal = {The Journal of Symbolic Logic},
volume = {78},
number = {4},
year = {2013},
pages = {1036-1054},
doi = {10.2178/jsl.7804020},
note = {Preprint arXiv:1012.2381}}

@book{BourbakiTopology1,author={Bourbaki, Nicolas},title={Elements of Mathematics. General Topology. Part 1},publisher={Hermann, Paris; Addison-Wesley, Reading, Mass.--London--Don Mills, Ont.},year={1966}}

@inproceedings{Bulatov2017,
  author        = {Bulatov, Andrei A.},
  title         = {A dichotomy theorem for nonuniform {CSPs}},
  booktitle     = {58th {IEEE} Annual Symposium on Foundations of
                   Computer Science ({FOCS})},
  pages         = {319--330},
  year          = {2017},
  eprint        = {1703.03021},
  archivePrefix = {arXiv},
  url           = {https://arxiv.org/abs/1703.03021}
}

@article{BJK2005,
  author  = {Bulatov, Andrei and Jeavons, Peter and Krokhin, Andrei},
  title   = {Classifying the complexity of constraints using finite algebras},
  journal = {SIAM Journal on Computing},
  volume  = {34},
  number  = {3},
  pages   = {720--742},
  year    = {2005},
  doi     = {10.1137/S0097539700376676}
}

@article{FV1998,
  author  = {Feder, Tom{\'a}s and Vardi, Moshe Y.},
  title   = {The computational structure of monotone monadic {SNP} and
             constraint satisfaction: A study through {Datalog} and group theory},
  journal = {SIAM Journal on Computing},
  volume  = {28},
  number  = {1},
  pages   = {57--104},
  year    = {1998},
  doi     = {10.1137/S0097539794266766}
}

@inproceedings{FP2026,
  author    = {Feller, Roman and Pinsker, Michael},
  title     = {Decidability of interpretability},
  booktitle = {41st Annual Symposium on Logic in Computer Science ({LICS} 2026)},
  series    = {Leibniz International Proceedings in Informatics ({LIPIcs})},
  volume    = {380},
  pages     = {42:1--42:26},
  publisher = {Schloss Dagstuhl -- Leibniz-Zentrum f{\"u}r Informatik},
  year      = {2026},
  doi       = {10.4230/LIPIcs.LICS.2026.42},
  eprint    = {2602.02302},
  archivePrefix = {arXiv}
}

@book{Hodges1997,author={Hodges, Wilfrid},title={A shorter model theory},publisher={Cambridge University Press},address={Cambridge},year={1997}}

@article{KPT05,author={Kechris, Alexander S. and Pestov, Vladimir G. and Todor{\v{c}}evi{\'c}, Stevo},title={Fra\"{i}ss\'e limits, {R}amsey theory, and topological dynamics of automorphism groups},journal={Geometric and Functional Analysis},volume={15},number={1},pages={106--189},year={2005}}

@misc{MarimonPinsker2025,author={Marimon, Paolo and Pinsker, Michael},title={A guide to topological reconstruction on endomorphism monoids and polymorphism clones},year={2025},eprint={2512.01086},archivePrefix={arXiv},note={Version of July 13, 2026}}

@article{MottetPinsker-cores-ramsey,
author = {Antoine Mottet and Michael Pinsker},
title = {Cores over {R}amsey structures},
journal = {The Journal of Symbolic Logic},
volume = {86},
number = {1},
year = {2021},
pages = {352-361},
doi = {10.1017/jsl.2021.6},
note = {Preprint arXiv:2004.05936}}

@article{NguyenVanThe-KPT-correspondence,
author = {Lionel {Nguyen Van Th\'{e}}},
title = {A survey on structural {R}amsey theory and topological dynamics with the {K}echris-{P}estov-{T}odorcevic correspondence in mind},
journal = {Zbornik Radova},
volume = {17},
number = {25},
year = {2015},
pages = {189-207},
note = {Selected Topics in Combinatorial Analysis, Mat. Inst. SANU, Beograd. Preprint arXiv:1412.3254}}

@article{Tsankov2012,
  author  = {Tsankov, Todor},
  title   = {Unitary representations of oligomorphic groups},
  journal = {Geometric and Functional Analysis},
  volume  = {22},
  pages   = {528--555},
  year    = {2012},
  doi     = {10.1007/s00039-012-0156-9}
}

@inproceedings{Zhuk2017,
  author        = {Zhuk, Dmitriy},
  title         = {A proof of {CSP} dichotomy conjecture},
  booktitle     = {58th {IEEE} Annual Symposium on Foundations of
                   Computer Science ({FOCS})},
  pages         = {331--342},
  year          = {2017},
  eprint        = {1704.01914},
  archivePrefix = {arXiv},
  url           = {https://arxiv.org/abs/1704.01914}
}

\end{document}